

\baselineskip=14pt
\parskip=10pt

\font\eightrm=cmr8 

\magnification=\magstephalf
\def\C{{\cal C}}
\def\N{{\cal N}}
\def\P{{\cal P}}

\def\M{{\cal M}}
\def\1{{\overline{1}}}
\def\2{{\overline{2}}}
\parindent=0pt
\overfullrule=0in

\def\frac#1#2{{#1 \over #2}}
\centerline
{\bf The ``Monkey Typing Shakespeare" Problem for Compositions }
\bigskip
\centerline
{\it Shalosh B. EKHAD and Doron ZEILBERGER}
\bigskip

{\bf Pattern Avoidance}

The {\it theme} of {\it pattern-avoidance}, i.e. the question 

{\it `How Many objects of such and such kind avoid such and such patterns?'}

is everywhere dense in combinatorics. Often the problems are so hard that the {\it mere} existence is highly non-trivial, for example
in Ramsey theory and Roth and Szemer\'edi's theorems.

{\eightrm [We understand the word `pattern' in a general sense, the notion of `permutation pattern' ([Wi],[HM]) is only one case.]}

One of the simplest such problems is when the objects are {\it words} in a fixed finite alphabet (say our very own A-Z, consisting
of $26$ letters) and the pattern to avoid is a fixed {\it consecutive} subword (say the collected works of the Bard arranged in
lexicographic order of titles, ignoring spaces and punctuation marks). It is easy to see that the
probability that a monkey typing $n$ letters randomly will {\it avoid} a run of the (long but finite) string consisting of
Willie's opus is less than $ \alpha^n$ for some fixed constant $\alpha$ {\it strictly} less than $1$, hence if
the monkey lives long enough he will sooner or later succeed. Of course it is latttttttttttttttttttttttttttttttter rather
than sooner, but we mathematicians pretend that we, and our monkeys, are immortal.

The problem of finding the {\it exact} probability of a random word in a fixed (finite) alphabet
avoiding a pre-determined string (or set of strings) as {\it consecutive} subwords 
(or equivalently counting the number of such words)
is handled by the {\it Goulden-Jackson Cluster Method} [GJ], beautifully exposited, and extended, in [NZ].

Another, much harder, example of the notion of {\it pattern-avoidance} is the one  used in
the very active field of {\it permutation patterns} [Wi], that has been beautifully extended
to compositions and words in Silvia Heurbach and Toufik Mansour's magnum opus [MH].

{\bf This article}

In this article we will also talk about compositions, but the notion of {\it pattern} is different.
If we define a pattern just as the {\it literal} occurrence of a consecutive subword, then we get
a straightforward generalization of the Goulden-Jackson set-up ([GJ] [NZ]).
In this article, we fill a much needed gap (Oops, we meant `fill a gap that was in great need of being filled')
and describe an algorithm, fully implemented in Maple, to enumerate  the problem of avoiding
(consecutive) {\it containment} of {\it compositions}, in a sense to be explained shortly.

The algorithm is fully implemented in the  Maple package {\tt Compositions.txt}, written by the second-named author, available from the url

{\tt http://sites.math.rutgers.edu/\~{}zeilberg/mamamrim/mamarimhtml/kof.html} \quad ,

where one can also find  extensive sample output files, created by the first-named author.

First we must recall the definition of {\it composition} and
then what we mean by a composition {\it avoiding} another composition by {\it containment}.

{\bf Compositions}

Recall that a {\it composition} of a non-negative integer $n$ is an {\bf ordered} list
$a_1\,a_2\, \dots\, a_k$ of {\bf positive} integers that add-up to $n$.

For example, the set of compositions of $4$ is

$$
\{ 4,13,31, 22, 112, 121, 211, 1111 \} \quad .
$$

The number of compositions of $n$ is famously $2^{n-1}$ for $n>0$ and $1$ for $n=0$. A quick way to see it
is by assigning the composition $a_1 a_2 \dots a_k$ to the subset of $\{1, ..., n-1\}$ given by
$\{a_1, a_1+a_2, \dots, a_1+ \dots + a_{k-1}\}$.
This is a bijection. Another way to prove this trivial fact is  via generating functions.
The generating function of a single entry is $\sum_{i=1}^{\infty} x^i = \frac{x}{x-1}$. The generating function
of compositions into exactly $k$ parts is, by `independence',  $( \frac{x}{x-1})^k$, hence the generating function of
all compositions is $\sum_{k=0}^{\infty} ( \frac{x}{x-1})^k=1/(1-  \frac{x}{x-1})=\frac{1-x}{1-2x}=1+ \frac{x}{1-2x}$.

{\bf When does a composition include another composition?}

{\bf Definition}: The composition $a_1 \dots a_k$ {\it includes} the composition $b_1 \dots b_s$ if
$k \geq s$ and  there exists an $i$, $1 \leq i \leq k-s+1$ such that
$$
b_1 \leq a_i \quad , \quad b_2 \leq a_{i+1}  \quad , \quad \dots  \quad , \quad b_s \leq a_{i+s-1} \quad .
$$

For example, the composition $222$ includes all the members of the set of compositions $\{ 1, 2, 11, 12, 21, 22\}$.

{\bf Definition}: The composition $a_1 \dots a_k$ {\it avoids} the composition $b_1 \dots b_s$ if
it does {\bf not} include  it. 

For example, the set of compositions of $4$ that do not include the composition $12$ is
$$
\{ 4,  31, 211, 1111 \} .
$$

{\bf Definition}: The composition $a=a_1 \dots a_k$ {\it avoids} the 
set of compositions $A$, if  it avoids every member of $A$. 

For example, 
the set of compositions of $4$ that do not include the members of the set $\{12 , 21  \}$  is
$$
\{ 4, 1111 \} .
$$

The {\bf million dollar question} that we will answer in this article is

{\bf Question}: Design an algorithm that inputs an arbitrary finite set of compositions, $A$,  {\it all of the same length},
and a variable $x$ and outputs  an explicit expression (a rational function of $x$) for
$$
\sum_{n=0}^{\infty} a(n) x^n \quad,
$$
where $a(n)$ is the number of compositions of $n$ that do {\bf not contain }  any of the members of $A$.

For the sake of simplicity, we only treat the case where all the members of the forbidden set are of the same length.
The extension to the more general case is left to the interested reader.

Note that when the set $A$ consists of a single composition of length $1$, i.e. $A=\{k\}$, then
the answer is very easy. It is the set of compositions where the parts are from $\{1, \dots, k-1\}$ and hence
the required generating function is simply $1/(1-x-x^2 - \dots -x^{k-1})$. For example when $A=\{3\}$ then
$a(n)=F_{n+1}$.

Before going on, let's give two sample theorems. Using our Maple package, the reader can generate
(potentially) infinitely many such theorems.

The first theorem enumerates compositions that never contain the composition $34543$.

It was obtained, {\it in less than a second}, by typing

{\tt InfoV( $\{$[3,4,5,4,3] $\}$,x,30): }

in the Maple package {\tt Compositions.txt }.

{\bf Theorem 1}:  Let $a(n)$ be the number of compositions of $n$ avoiding, as a subcomposition, $34543$, then
$$
\sum_{n=0}^{\infty} a(n)\, x^n \, = \,
$$
$$
-{\frac {1-4\,x+6\,{x}^{2}-4\,{x}^{3}+{x}^{4}+{x}^{16}+{x}^{13}-{x}^{14}+{x}^{9}-2\,{x}^{10}+{x}^{11}+{x}^{5}-3\,{x}^{6}+3\,{x}^{7}-{x}^{8}}{{x}^{18}+{x}^{
17}-{x}^{16}+2\,{x}^{14}-{x}^{13}-2\,{x}^{11}+3\,{x}^{10}-{x}^{9}+2\,{x}^{8}-5\,{x}^{7}+4\,{x}^{6}-{x}^{5}-2\,{x}^{4}+7\,{x}^{3}-9\,{x}^{2}+5\,x-1}} .
$$

The first $31$ terms of $a(n)$,  starting at $n=0$, are:
$$
1, 1, 2, 4, 8, 16, 32, 64, 128, 256, 512, 1024, 2048, 4096, 8192, 16384, 32768, 65536, 131072, 262143, 
$$
$$
524281, 1048546, 2097050, 4194001, 8387784, 16775108,  33549270, 67096623, 134189393, 268371074, 536726740 .
$$
The limit of $a(n+1)/a(n)$ as $n$ goes to infinity is $1.99994300442\dots$, and $a(n)$ is asymptotic to
$ (0.50029301491\dots)\, (1.99994300442\dots)^n$.

The second theorem enumerates compositions that do not contain any of the members of  $\{ 252, 343, 424 \}$.
It was obtained by typing

{\tt InfoV( $\{$[2,5,2],[3,4,3],[4,2,4]  $\}$,x,30): }

in the Maple package {\tt Compositions.txt }.

{\bf Theorem 2}:  Let $a(n)$ be the number of compositions of $n$ avoiding, as subcompositions, all the three members the set
$\{ 252, 343, 424 \}$, then
$$
\sum_{n=0}^{\infty} a(n)\, x^n \, = \,
$$
$$
-{\frac {{x}^{17}+3\,{x}^{14}+{x}^{13}-3\,{x}^{11}+{x}^{10}+{x}^{6}-{x}^{5}+{x}^{4}+{x}^{2}-2\,x+1}{{x}^{18}+3\,{x}^{15}+2\,{x}^{14}-2\,{x}^{13}-2\,{x}^{12}
+3\,{x}^{11}-3\,{x}^{10}+{x}^{8}+{x}^{7}-{x}^{6}+2\,{x}^{5}-{x}^{4}-2\,{x}^{2}+3\,x-1}}.
$$

The first $31$ terms of $a(n)$,  starting at $n=0$, are:
$$
1, 1, 2, 4, 8, 16, 32, 64, 128, 255, 505, 998, 1971, 3893, 7697, 15223, 30113, 59575, 117861, 233164, 
$$
$$
461250, 912423, 1804882, 3570257, 7062369, 13970211, 27634848,    54665348, 108135332, 213906125, 423134791 .
$$
The limit of $a(n+1)/a(n)$ as $n$ goes to infinity is $1.9781317474\dots$, and $a(n)$ is asymptotic to
$(.54805291269\dots)\,(1.9781317474\dots)^n$.

{\bf The Cluster Method for Compositions}

Our algorithm is an adaptation of the Goulden-Jackson cluster method as described in [NZ].

We will use generating functions. The weight of a composition, $C$, is defined by
$$
Weight(C):=  \, x^{Sum(C)} \quad .
$$
The weight-enumerator of the set of compositions, as noted above, is $1/(1-(\frac{x}{1-x}))=\frac{1-x}{1-2x}$.
Let's reprove this trivial fact, in order to motivate our algorithm.

Any composition is either empty (weight $1$) or starts with a {\bf positive integer}. Removing the first
entry results in a brand-new composition, i.e.
$$
C \, = \, EmptyComposition  \quad or \quad C\, = \, i\, C' \quad,
$$
for some positive integer $i$, and  some smaller (possibly empty) composition $C'$. 
Let $\C$ be the set of compositions, and $\P=\{1,2,3,\dots\}$ be the
set of {\bf strictly positive} integers. We have
$$
\C = \{EmptyComposition\} \, \bigcup \, \P \times \C \quad .
$$
Taking weights on both sides (the weight of a set is the sum of the weights of its members), we have
$$
Weight(\C) = Weight(\{EmptyComposition\}) + Weight(\P)\, Weight(\C)  \quad .
$$
Since  $Weight(\{EmptyComposition\})=1$ and  $Weight(\P)=\sum_{i=1}^{\infty} x^i = \frac{x}{1-x}$, we have
the {\it algebraic} equation
$$
Weight(\C) = 1 + \frac{x}{1-x} \, Weight(\C)  \quad ,
$$
and solving  for $Weight(\C)$ we get that it indeed equals $1/(1- \frac{x}{1-x})$.

Our problem is as follows. We are given a finite set of compositions, of the {\bf same length } (that is our simplifying assumption), 
let's call it $A$, and we would like to find
the weight-enumerator of the set of compositions, let's call it $\C(A)$, that do {\bf not} contain any of the members of
$A$. For simplicity of exposition, let's assume, until further notice, that the set of offenders is a singleton.
Later we will describe how to extend it to the case when we want to avoid, by containment,  more than one fixed composition.

In other words, we would like to find a way of computing
$$
Weight(\C(A))(x) \quad .
$$
Being the weight-enumerator of the set $\C(A)$, extracting the coefficient of $x^n$ from its Maclaurin expansion would
give us the exact number of compositions of $n$ that do {\bf not} contain any of the members of $A$.

Rather than counting {\it good guys}, that is very hard, we will do {\it signed-weighted counting} of
pairs $(C,S)$, where $C$ is {\it any} composition, and $S$ is a subset of its set of `crimes', i.e. a
subset of its set of offending containments, that in turn, is a certain subest of $A$ (possibly empty, if $C$ is a saint, and
possibly the whole of $A$, if $C$ is an arch-criminal).
Let's call this much larger set $\hat{\C}(A)$, and define
$$
\overline{Weight}(C,S):= x^{Sum(C)}(-1)^{|S|} \quad .
$$

Using the following two extremely deep identities
$$
1+(-1) \, = \, 0 \quad ,
$$
and
$$
0^a \,= \, \cases{ 1 \quad if \quad a=0 ; \cr
                  0 \quad if \quad a>0 .} \quad ,
$$
it follows, just as in the case of Goulden-Jackson, that
$$
Weight(\C(A)) \, = \, \overline{Weight}(\hat{\C}(A)) \quad .
$$
We are left with the task of computing the right side,   $\overline{Weight}(\hat{\C}(A))$.

let's illustrate the method by using a concrete example, where the set of offenders, $A$, 
consists of the single composition $232$. In other words $A=\{232\}$.
How do we compute  $\overline{Weight}(\hat{\C}(\{232\}))$? (and hence
$Weight(\C(\{232\}))$?

Consider the composition $14542351$, whose $Weight$  is $x^{1+4+5+4+2+3+5+1}=x^{25}$. It has a total of three `crimes'
$$
\matrix{ 1 & 4 & 5 & 4 & 2 & 3 & 5 & 1\cr
         \_ & \_ & \_& \_ & \_ & \_ & \_ & \_  \cr 
           & 2 & 3 & 2 &   &   &   &  }
\quad \quad, \quad \quad
\matrix{ 1 & 4 & 5 & 4 & 2 & 3 & 5 & 1\cr
         \_ & \_ & \_& \_ & \_ & \_ & \_ & \_  \cr
           &  &  2 & 3 & 2  &   &   &  }
\quad \quad, \quad \quad
\matrix{ 1 & 4 & 5 & 4 & 2 &  3 & 5 & 1\cr
         \_ & \_ & \_& \_ & \_ & \_ & \_ & \_  \cr
           &   &   &   & 2  & 3 & 2  &  } \quad .
$$

These give rise to $2^3$ members of $\hat{\C}(\{232\})$. For example, when  $S=\emptyset$:
$$
\matrix{ 1 & 4 & 5 & 4 & 2 & 3 & 5 & 1\cr
         \_ & \_ & \_& \_ & \_ & \_ & \_ & \_  \cr
           &  &  &  &   &   &   &  } \quad,
$$
whose weight is $(-1)^0\, x^{25}=x^{25}$. 
Another example, is when $S$ consists of the violations that start at the 2nd and 3rd entries of the underlying composition $C=14542351$:
$$
\matrix{ 1 & 4 & 5 & 4 & 2 & 3 & 5 & 1\cr
         \_ & \_ & \_& \_ & \_ & \_ & \_ & \_  \cr
           & 2 & 3 & 2 &   &   &   &   \cr
           &   & 2 & 3 & 2 &   &   &
}\quad,
$$
whose weight is $(-1)^2 x^{25}=x^{25}$. Yet another example, out of the eight possibilities , is
when $S$ consists of violations starting at the $2$nd and $5$th entries of $14542351$:
$$
\matrix{ 1 & 4 & 5 & 4 & 2 & 3 & 5  & 1\cr
         \_ & \_ & \_& \_ & \_ & \_ & \_ & \_  \cr
           & 2 & 3 & 2 &   &   &    &   \cr
           &   &   &   & 2 & 3 & 2   &
}\quad,
$$
whose weight is also $(-1)^2\, x^{25}=x^{25}$. Finally, don't forget the case where $S$ is the full set of `crimes':
$$
\matrix{ 1 & 4 & 5 & 4 & 2 & 3 & 5 & 1\cr
         \_ & \_ & \_& \_ & \_ & \_ & \_ & \_  \cr
           & 2 & 3 &  2 &   &   &   &   \cr
          &    & 2  & 3 & 2 &     &    & \cr
           &   &  &     & 2 &   3  & 2   &
}\quad,
$$
whose weight is $(-1)^3\, x^{25}=- x^{25}$.

Let's analyze the anatomy of a typical member of $\hat{\C}(A)$, $(C,S)$.
We have one of the following three cases

$\bullet$  It is the element $(EmptyComposition, \emptyset)$.

$\bullet$ The first entry of $C$ does not belong to any offenders. This case is isomorphic to $\P \times \hat{\C}(A)$

$\bullet$ The first entry belongs to an offender. This offender may not overlap with any other offenders that
start later on, or may. Continuing to examine whether the last offender in the chain overlaps with yet another
offender, eventually we will stop, getting a `{\bf pre-cluster}', after which we have a brand-new  (shorter) member of
$\hat{\C}(A)$. Let's denote by $\M(A)$ this set of pre-clusters.
It is the subset of $\hat{\C}(A)$ where every entry of the underlying composition
belongs to at least one offender, and every offender overlaps with at least another offender.

Still using $A=\{232\}$ as our running example, here is a member of $\M(A)$:

$$
\matrix{ 5 & 4 & 5 & 4 &  5 & 3  & 5 & 6 & 4\cr
         \_ & \_ & \_& \_ & \_ & \_ & \_ & \_   & \_  \cr 
         2 & 3 & 2 &   &    &    &   &   &    \cr
           &   & 2 & 3 &  2 &    &   &   &    \cr
           &   &   & 2 &  3 &  2 &   &   &    \cr
           &   &   &   &    & 2  & 3 &  2&   \cr
           &   &   &   &    &    & 2 &  3& 2 
}\quad \quad .
$$
On the other hand
$$
\matrix{ 5 & 4 & 5 & 4 &  5 & 3  & 5 & 6 & 4\cr
         \_ & \_ & \_& \_ & \_ & \_ & \_ & \_   & \_  \cr 
         2 & 3 & 2 &   &    &    &   &   &    \cr
           &   &   & 2 &  3 & 2   &   &   &   } \quad ,
$$
is {\bf not} a member of $\M(\{232\})$, since the second offending $232$ does not overlap the first one, and hence 
can be decomposed to
$$
\left [ \matrix{ 5 & 4 & 5 \cr
         \_ & \_ & \_ \cr
         2 & 3 & 2  \cr} \right ]
\,
\left [ \matrix{  4 &  5 &  3  & 5 & 6 & 4\cr
        \_ & \_ & \_ & \_ & \_   & \_  \cr 
          2 &  3 &  2  &   &   &   } \right ]\quad .
$$

We get the following {\it grammar} for our set $\hat{\C}(A)$.
$$
\hat{\C}(A) = \, \, \{(EmptyComposition, \emptyset)\} \,\, \bigcup \,\, \P \times \hat{\C}(A) \,\,  \bigcup \,\,  \M(A) \times  \hat{\C}(A) \quad .
$$
Applying $\overline{Weight}$, and abbreviating $\overline{Weight}(\hat{\C}(A))$ (our object of desire) to $F(x)$,
we have
$$
F(x) \, = \, 1 + \frac{x}{1-x} \, F(x) + \overline{Weight}(\M(A)) F(x) \quad .
$$
Solving for $F(x)$, and abbreviating  $\overline{Weight}(\M(A))$ to $G(x)$, 
we get
$$
F(x)= \frac{1}{1 \, - \, \frac{x}{1-x} \, - \, G(x)} \quad .
$$

It remains to compute the  $\overline{Weight}$-enumerator of the set of pre-clusters, $\M(A)$, alias $G(x)$.

Let's forget for a second the underlying composition (the top row in the above examples), and focus on the offenders.
Let's call this pre-cluster with the top row (underlying composition) removed a {\bf cluster}.
We will denote the set of clusters by $\N(A)$. 

For example, still using $A=\{232\}$, here is an example of a cluster:
$$
\matrix{  2 & 3 & 2 &   &    &    &   &   &    \cr
           &   & 2 & 3 &  2 &    &   &   &    \cr
           &   &   & 2 &  3 &  2 &   &   &    \cr
           &   &   &   &    &  2  & 3 &  2 &   \cr
           &   &   &   &    &     & 2 &  3 & 2 
}\quad \quad .
$$
Let's call  the total number of columns the {\it width}. In the above example, the width is $9$.
What compositions may serve as the underlying composition (top row) of such a cluster?
Obviously, its number of columns must be equal to the width, and its respective entries, must be
$\geq$ to the maximum of the entries of the offenders in the corresponding columns, where
we ignore the empties (or replace them by $0$, note that every column has at least one non-empty
entry, or else it would not be a cluster). In the above example, the maxima of the nine columns are
$$
2 \, 3 \, 2 \, 3 \, 3 \, 2 \, 3 \, 3 \, 2 \quad .
$$

We will call this the {\it Skyline} of the cluster.

For  the sake of convenience, let's represent a cluster as  matrix, with empties replaced by $0$, then
the above cluster is written
$$
\matrix{  2 & 3 & 2 &  0 &  0  & 0   & 0  & 0  & 0   \cr
          0 & 0  & 2 & 3 &  2 &  0  &  0 &  0  & 0   \cr
          0 & 0  & 0  & 2 &  3 &  2 &  0  & 0  & 0   \cr
          0 & 0  & 0  & 0  & 0   &  2  & 3 &  2 &0   \cr
          0 & 0  & 0  & 0  & 0   &  0   & 2 &  3 & 2 
}\quad \quad .
$$
In terms of this associated matrix (where empty spaces are replaced by $0$), we 
get  a certain $k \times r$ matrix ($k$ is the number of offenders, and $r$ is the width),
the Skyline is defined as follows.
$$
 Skyline(L)_j:= max(\{ L_{ij} \, | \, 1 \leq i \leq k \} ) \quad , \quad 1 \leq j \leq r .
$$

Note that the underlying compositions that can be put on top of a cluster $L$  must satisfy $C_j \geq Skyline(L)_j$ for all $1\leq j \leq r$.
The $\overline{Weight}$-enumerator of these is simply $x^{Sum(Skyline(L))} \, (\frac{1}{1-x})^r (-1)^k$.

This leads us to define yet-another weight, this time on clusters, introducing an auxiliary variable $t$:
$$
Poids(L)(x,t) :=x^{Sum(Skyline(L))} \, (-1)^k \, t^r \quad .
$$
We have
$$
G(x)= \overline{Weight}(\M(A))= Poids(\N(A)) (x, \frac{1}{1-x}) \quad .
$$

We are left with the task of computing  $Poids(\N(A))(x,t)$.

A natural approach would be to  break the set $\N(A)$ into {\it states}, and relate them to each other.
Continuing with the assumption that our set of offenders is a singleton, $\{A\}$, where $A$  is  of length $a$, say, ($A=232$ in our example, so $a=3$),
the state of a cluster is the list consisting of the first $a$ entries of its Skyline.
For example the state of the cluster
$$
\matrix{2 & 3 & 2 & 0 \cr
        0 & 2 & 3 & 2} \quad ,
$$     
is $233$, while the state of the cluster
$$
\matrix{2 & 3 & 2 & 0 & 0 \cr
        0 & 0 & 2 & 3 & 2} \quad ,
$$     
is $232$. It is readily seen that in this case these are the only states. In general the computer can easily
determine the set of states, by generating all the clusters of width $\leq 2a-1$ (recall that $a$ is the width of $A$), and
extracting the first $a$ entries.

For each state, $s$, let $B_s(x,t)$ be the Poids-enumerator of the set of clusters that belong to state $s$.
Given such a cluster, if it only has one row, then it must be $A$, and its Poids is $(-1)x^{Sum(A)} t^a$.
If it has more than one row, we look at all the possibilities that it can overlap with the second row
(i.e. the number of starting zeros (empties) in the second row, and their state).
These are the `children states'. In fact it is much easier, to go the other way. For any state
look at the `parent states', i.e. the one obtained by putting $A$ at the top row in all the $a-1$ possible ways.

For example, for the state $232$, we can have
$$
\matrix{2 & 3 & 2 & 0 \cr
\_ & \_ & \_ & \_ \cr
 0 & 2 & 3 & 2
} \quad ,
$$
whose state is $233$, and
$$
\matrix{2 & 3 & 2 & 0 & 0\cr
\_ & \_ & \_ & \_ & \_ \cr
        0 & 0 & 2 & 3 & 2} \quad ,
$$
whose state is $232$. Hence the parents of the state $232$ are $233$ and $232$.

On the other hands for the  state $233$, we can have
$$
\matrix{2 & 3 & 2 & 0 \cr
\_ & \_ & \_ & \_ \cr
        0 & 2 & 3 & 3} \quad ,
$$
whose state is $233$, and
$$
\matrix{2 & 3 & 2 & 0 & 0\cr
\_ & \_ & \_ & \_ & \_ \cr
        0 & 0 & 2 & 3 & 3}
$$
whose state is $232$. Hence the parents of the state $233$ are also  $233$ and $232$.

Recall that $B_s(x,t)$ is the Poids-enumerator of the clusters with state $s$.
For each state, we need to set-up an equation.

Let's first find an equation for $B_{232}$. The children of the state $232$ are $232$ and $233$. 
If the child is $232$ then it {\bf must} be as follows
$$
\matrix{2 & 3 & 2 & 0 & 0 \cr
        \_ & \_ & \_ & \_ & \_ \cr
        0 & 0 & 2 &  3& 2} \quad .
$$
Removing the first row costs $2+3=5$ units to the Skyline, and shrinks the width by $2$, so the contribution 
to $B_{232}$ from this scenario (where the child is $232$) is $x^5 t^2 (-1)B_{232}$.

If the child is $233$ then it {\bf must} be as follows
$$
\matrix{2 & 3 & 2 & 0 & 0 \cr
        \_ & \_ & \_ & \_ & \_ \cr
        0 & 0 & 2 &  3& 3} \quad .
$$
Removing the first row costs $2+3=5$ units to the Skyline, and shrinks the width by $2$, so the contribution 
to $B_{232}$ from this scenario (where the child is $233$) is $x^5 t^2 (-1)B_{233}$.

In addition it may be just a one-row cluster, whose Poids is $-x^{2+3+2} \, t^3= -x^7\,t^3$.

Hence the equation for $B_{232}$ is
$$
B_{232} \, = \, -x^7\, t^3 -x^5\,t^2\,B_{232}  -x^5\,t^2\,B_{233}  \quad .
$$

Next, let's find an equation for $B_{233}$. The children of the state $233$ are also $232$ and $233$. 
If the child is $232$ then it {\bf must} be as follows
$$
\matrix{2 & 3 & 2 & 0  \cr
\_ & \_ & \_ & \_  \cr
        0 & 2 & 3 & 2} \quad .
$$
Removing the first row costs $2+3-2=3$ units to the Skyline, and shrinks the width by $1$, so the contribution 
to $B_{233}$ from this scenario (where the child is $232$) is $x^3 t (-1)B_{232}$.

If the child is $233$ then it {\bf must} be as follows
$$
\matrix{2 & 3 & 2 & 0  \cr
       \_ & \_ & \_ & \_  \cr
        0 & 2 & 3 &  3}
$$
removing the first row costs $3$ units to the Skyline, and shrinks the width by $1$, so the contribution 
to $B_{232}$ from this scenario (where the child is $233$) is $x^3 t (-1)\,B_{233}$.

Since $233$ is not a violation, we don't have to add its Poids.

Hence the equation for $B_{233}$ is
$$
B_{233} \, = \, -x^3\,t \,B_{232}  -x^3\,t \,B_{233}  \quad .
$$

We have to solve the following system of two linear equations with two unknowns $\{ B_{232}, B_{233} \}$:
$$
\{ \,
B_{232} \, = \, -x^7\, t^3 -x^5\,t^2\,B_{232}  -x^5\,t^2\,B_{233}  \quad , \quad
B_{233} \, = \, -x^3\,t \,B_{232}  -x^3\,t \,B_{233}  \, \} \quad .
$$
Solving them gives
$$
B_{232} \, = \, -{\frac { \left( 1+t{x}^{3} \right) {t}^{3}{x}^{7}}{1+t{x}^{3}+{t}^{2}{x}^{5}}} \quad , \quad
B_{233} \, = \,{\frac {{t}^{4}{x}^{10}}{1+t{x}^{3}+{t}^{2}{x}^{5}}} \quad .
$$

Hence
$$
Poids(\N(\{232\}))(x,t)=B_{232}+B_{233} = -{\frac {{t}^{3}{x}^{7}}{1+t{x}^{3}+{t}^{2}{x}^{5}}}  \quad .
$$
Hence
$$
G(x)=\overline{Weight}(\M(\{232\})= Poids(\N(\{232\}))(x,\frac{1}{1-x})
= {\frac {{x}^{7}}{ \left( 1-2\,x+{x}^{2}+{x}^{3}-{x}^{4}+{x}^{5} \right)  \left( -1+x \right) }} \quad .
$$
Finally, our object of desire, $F(x)=\overline{Weight}(\hat{\C}(A))(x)$, alias $Weight(\C(A))$, is
$$
F(x)= \frac{1}{1-\frac{x}{1-x} - G(x)} \, = \,
      -{\frac {1-2\,x+{x}^{2}+{x}^{3}-{x}^{4}+{x}^{5}}{{x}^{6}-{x}^{5}+2\,{x}^{4}-{x}^{3}-2\,{x}^{2}+3\,x-1}} \quad .
$$

{\bf Avoiding many compositions (by containment)}

The above method can be easily modified to handle multiple  offenders. The clusters now can have several
violations starting at any location, and the effect of the Skyline on a given state is determined by
looking at the Skyline of the subset of violations starting at the same column. The readers are
welcome to look at the source code of procedure {\tt GFset(S,x)} in the Maple package
{\tt Compositions.txt}, where this case is implemented. So far we only handle sets of offenders all of the same length.
Our readers are welcome to extend it to the more general case where the offending compositions may be of different lengths.

{\bf Keeping track of the number of occurrences}

So far we described how to find the exact number, let's call it $a_S(n)$ of compositions of length $n$, that {\bf avoid}
(i.e. do not contain) any members of the set of compositions $S=\{ C_1, ..., C_r \}$ . We gave an efficient
algorithm, implemented in the Maple package {\tt Compositions.txt}, to explicitly find the generating function,
let's call $f_S(x)$
$$
f_S(x) \, := \, \sum_{n=0}^{\infty} \, a(n) \, x^n \quad .
$$

Note that, out of laziness, (so far) we only implemented the case where all the members of the offending set, 
$S=\{ C_1, ..., C_r \}$, are of the same length.

But it is very hard to stay out of trouble. Suppose that you want to find the exact number, let's call it
$A_S(n; c_1, \dots, c_r)$ of compositions of $n$ that contain 

$C_1$, $c_1$ times,  $C_2$, $c_2$ times,  \dots, $C_r$, $c_r$ times.

Of course, our former quantity, $a_S(n)$ is just the special case $c_1=0, c_2=0 , \dots, c_r=0$, i.e.
$$
a_S (n) \, = \,  A_S(n; 0,0, \dots, 0) \quad .
$$

All the information about the discrete function $A_S(n; c_1, \dots, c_r)$, with $1+r$ {\it discrete} variables, is
{\bf encapsulated} in the multi-variable rational function, with $1+r$ `{\it continuous}' variables
$x,X_1, \dots, X_r$
$$
F_S(x; X_1, \dots, X_r) \, := \, \sum_{n=0}^{\infty}  \sum_{c_1=0}^{\infty}  \dots \sum_{c_r=0}^{\infty} 
\, A(n; c_1, \dots, c_r)\, x^n {X_1}^{c_1} \cdots {X_r}^{c_r}  \quad .
$$

Of course $f_S(x)=F_S(x; 0,0 , \dots, 0)$.

The beauty of the cluster method is that a very tiny tweak in the former algorithm yields a way to compute $F_S$.
Rather then use the deep identity $0=1+(-1)$ we use the only slightly deeper identity
$X=1+(X-1)$ applied to $X=X_1$, $X=X_2$, $\dots$, $X=X_r$. This entails  redefining the {\it Poids} of a cluster $C$, to be

$x^{Skyline(C)}$ (as before) 

{\bf multiplied by} $(X_1-1)^{NumberOf\,C_1}$ \quad ,

{\bf multiplied by} $(X_2-1)^{NumberOf\,C_2}$ \quad ,

$\dots$ \quad ,

{\bf multiplied by} $(X_r-1)^{NumberOf\,C_r}$ \quad .

The set of equations is modified accordingly, and Maple solves it. Of course, now it
takes much longer, since we have so many more symbols, but the principle is the same.

This is implemented (still under the simplifying assumption of the members of the offending set all
of the same length) in the Maple package {\tt CompositionsPlus.txt}, also available from the front
of this article, where there are also sample input and output files.

Using the method of [Z], (that have been included in this package) we can do {\it statistical analysis}
of the random variables `Number of occurrences of $C_i$', for the various members  of $S$, defined
on the sample space of all compositions of $n$, as well as how they interact.

We do it by using the multi-variable generating function $F_S$.

By taking partial derivatives, and then setting all the variables $X_i$ to be $1$,
we find expressions for the expectation, variance, (these are always linear in $n$)
and higher moments (certain polynomials in $n$).
We can also find mixed moments for any set of such
random variables, in particular, the {\it asymptotic correlation}, and confirm that for any such pair, these
random variables are {\it joint asymptotically normal}, alas (of course) , {\bf not} independently so.
Using the asymptotic correlation one can confirm this by computing the mixed moments of the corresponding
bi-variate normal distribution with correlation $\rho$, 
$$
\frac{1}{2\pi \, \sqrt{1-\rho^2}} \, e^{-x^2/2 -y^2/2 + \rho\,x\,y} \quad .
$$
Our Maple package does that automatically (to any desired order).

Just to cite one example, typing

{\tt InfoX2V([2,3,4],[4,3,2],x,X,Y,n,6):}

in the Maple package {\tt CompositionsPlus.txt}, yields the following theorem.

{\bf Theorem 3}: The following statements are true.

$\bullet$ Let $a(n)$ be the number of compositions of $n$ that contain  neither $234$ nor $432$, then
$$
\sum_{n=0}^{\infty} \, a(n) \, x^n \, = \,
-{\frac {{x}^{16}+{x}^{15}+{x}^{12}+2\,{x}^{10}-{x}^{7}+{x}^{5}-{x}^{4}-{x}^{2}+2\,x-1}{{x}^{17}+{x}^{16}+{x}^{13}+2\,{x}^{11}-{x}^{10}+{x}^{9}-{x}^{8}+{x}^
{7}-2\,{x}^{5}+{x}^{4}+2\,{x}^{2}-3\,x+1}} \quad .
$$

$\bullet$ $a(n)$ is asymptotic to $(0.548269839581\dots) \cdot (1.976902834153\dots)^n$

$\bullet$ Let $A(n; c,d)$ be the number of compositions of $n$ that contain exactly $c$ occurrences of $234$ and
$d$ occurrences of $432$, then
$$
\sum_{n=0}^{\infty} \sum_{c=0}^{n-3} \sum_{d=0}^{n-3} \, A(n;c,d) \, x^n \, X^c \, Y^d \, = \, \frac{Numer(x,X,Y)}{Denom(x,X,Y)} \quad, where
$$
$$
Numer(x,X,Y) \, = \,
-1+{x}^{12}-2\,{x}^{11}+{x}^{3}+2\,{x}^{5}-{x}^{6}-{x}^{4}-{x}^{7}-3\,{x}^{2}+{x}^{8}+2\,{x}^{10}-{x}^{17}+{x}^{15}{X}^{2}-{x}^{13}{X}^{2}
$$
$$
-{x}^{12}X-{x}^{12
}Y-{x}^{13}{Y}^{2}+{x}^{15}{Y}^{2}-{x}^{17}{Y}^{2}-{x}^{17}{X}^{2}-2\,{x}^{15}X-2\,{x}^{15}Y+2\,{x}^{11}Y-2\,{x}^{10}X+2\,{x}^{13}Y+2\,{x}^{11}X
$$
$$
-2\,{x}^{10}
Y+2\,{x}^{13}X+2\,{x}^{17}Y+2\,{x}^{17}X+{x}^{15}+3\,x+{x}^{7}XY-{x}^{13}+{x}^{6}XY-{x}^{8}XY+{x}^{4}XY+{x}^{15}{X}^{2}{Y}^{2}-{x}^{17}{X}^{2}{Y}^{2}
$$
$$
+{x}^{13}X{Y}^{2}+{x}^{13}{X}^{2}Y+{x}^{12}XY-3\,{x}^{13}XY-2\,{x}^{11}XY+2\,{x}^{10}XY-2\,{x}^{15}{X}^{2}Y-2\,{x}^{15}X{Y}^{2}
$$
$$
+4\,{x}^{15}XY+2\,{x}^{17}{X}^{2}Y-2\,{x}^{5}XY-4\,{x}^{17}XY+2\,{x}^{17}X{Y}^{2} \quad ,
$$
and
$$
Denom(x,X,Y) \, = \,
-1+2\,{x}^{12}-3\,{x}^{11}+2\,{x}^{3}+3\,{x}^{5}-2\,{x}^{6}-{x}^{4}-{x}^{7}-5\,{x}^{2}+2\,{x}^{8}+2\,{x}^{10}+2\,{x}^{16}Y+2\,{x}^{16}X-2\,{x}^{18}Y
$$
$$
-2\,{x}^{14}X-2\,{x}^{14}Y-2\,{x}^{18}X-2\,{x}^{12}X-2\,{x}^{12}Y+{x}^{14}{X}^{2}+{x}^{9}X+{x}^{14}{Y}^{2}+{x}^{18}{Y}^{2}-{x}^{16}{Y}^{2}+{x}^{18}{X}^{2}
+{x}^{9}Y-{x}^{16}{X}^{2}
$$
$$
+3\,{x}^{11}Y-2\,{x}^{10}X+{x}^{13}Y+3\,{x}^{11}X-2\,{x}^{10}Y+{x}^{13}X+4\,x+2\,{x}^{16}X{Y}^{2}+{x}^{7}XY+4\,{x}^{18}XY-4\,{x}^{16}XY
$$
$$
-{x}^{13}-2\,{x}^{18}X{Y}^{2}+2\,{x}^{6}XY-2\,{x}^{18}{X}^{2}Y-2\,{x}^{8}XY+2\,{x}^{16}{X}^{2}Y+3\,{x}^{14}XY+{x}^{4}XY+2\,{x}^{12}XY-2\,{x}^{9}+{x}^{14}
$$
$$
-{x}^{13}XY-3\,{x}^{11}XY+2\,{x}^{10}XY-{x}^{14}{X}^{2}Y-{x}^{14}X{Y}^{2}-{x}^{16}{X}^{2}{Y}^{2}+{x}^{18}{X}^{2}{Y}^{2}-3\,{x}^{5}XY+{x}^{18}-{x}^{16} \quad .
$$

$\bullet$ The expectation and variance of the random variables `number of occurrences of $234$', and 
`number of occurrences of $432$', are both (obviously they are the same)
$$
\frac{n}{128} \quad , \quad \frac{147}{16384} \, n \, - \frac{1439}{16384}  \quad .
$$

$\bullet$ The asymptotic correlation is $\frac{71}{147} \, = \, 0.482993197\dots$, and the joint asymptotic normality (with that
correlation) is confirmed up to the sixth mixed moments (not that we had any doubts).

{\bf Encore: The asymptotic growth constants for all compositions up that of 6}

This is an excerpt from the output file 

{\tt http://sites.math.rutgers.edu/\~{}zeilberg/tokhniot/oCompositions4.txt}, 

ranking them according to the asymptotic growth constants of the sequences enumerating compositions that do not contain them.
We only list one of them in case of ties (due to trivial  equivalence).

$n=2$: $2\, (1)$  .

$n=3$: $12 \, (1)$, $3\, ( 1.6180339887498948482)$.

$n=4$: $112 \, (1)$, $13 \,(1.6180339887)$, $22 \, (1.7548776662)$, $4 \, ( 1.8392867552)$  .

$n=5$: $1112 \, (1)$, $113 \, (1.6180339887)$, $212 \, (1.7548776662)$, $14 \, (1.83928675521)$, $23 \, (1.86676039917)$, $5 \, ( 1.92756197548)$ .

$n=6$: $11112\, (1)$, $1113 \, (1.6180339887)$, $ 2112\, (1.7548776662)$, 
$114 \,(1.839286755214)$, $213 \,(1.866760399)$, $222 \,(1.908790738787)$, $15\, (1.92756197548)$,
$24 \, (1.93318498189952)$, $33 \, (1.9417130342786)$, $6\, (1.965948236645)$.

For the ranking for the compositions of up to $11$, see the above output file.

{\bf References}

[GJ] Ian Goulden and David Jackson, {\it An inversion theorem for cluster decomposition of sequences with distinguished 
subsequences},   J. London Math. Soc.(2) {\bf 20} (1979), 567-576. 

[HM] Silvia Heurbach and Toufik Mansour, {\it ``Combinatorics of Compositions and Words''},
Chapman and Hall/CRC, 2009.

[NZ] John Noonan and Doron Zeilberger, {\it The Goulden-Jackson cluster method: extensions, applications, and implementations},
 J. Difference Eq. Appl. {\bf 5} (1999), 355 -377.  \hfill\break
{\tt http://sites.math.rutgers.edu/\~{}zeilberg/mamarim/mamarimhtml/gj.html} \quad .

[Wi] Wikipedia, The Free Encyclopedia. {\it Permutation Pattern}, Retrieved 17:53, January 3, 2019 \hfill\break
{\tt https://en.wikipedia.org/wiki/Permutation\_pattern} \quad .

[Z] Doron Zeilberger,
{\it Automated derivation of limiting distributions of combinatorial random variables whose generating functions are rational},
The Personal Journal of Shalosh B. Ekhad and Doron Zeilberger, Dec. 24, 2016. \hfill\break
{\tt http://sites.math.rutgers.edu/\~{}zeilberg/mamarim/mamarimhtml/crv.html } \quad .

\bigskip
\hrule
\bigskip
Shalosh B. Ekhad, c/o D. Zeilberger, Department of Mathematics, Rutgers University (New Brunswick), Hill Center-Busch Campus, 110 Frelinghuysen
Rd., Piscataway, NJ 08854-8019, USA. \hfill\break
Email: {\tt ShaloshBEkhad at gmail dot com}   \quad .
\bigskip
Doron Zeilberger, Department of Mathematics, Rutgers University (New Brunswick), Hill Center-Busch Campus, 110 Frelinghuysen
Rd., Piscataway, NJ 08854-8019, USA. \hfill\break
Email: {\tt DoronZeil at gmail  dot com}   \quad .
\bigskip
\hrule
\bigskip
Exclusively published in the Personal Journal of Shalosh B.  Ekhad and Doron Zeilberger and arxiv.org \quad .
\bigskip
Written: Jan. 12, 2019.
\end